\documentclass[letterpaper, 10 pt, conference]{ieeeconf}  

\IEEEoverridecommandlockouts                              
\usepackage{cite,amsmath,amssymb}

\usepackage{tikz}
\usepackage{empheq}
\usepackage{algorithmic}
\usepackage{graphicx}
\usepackage{textcomp}
\usepackage{xcolor}
\usetikzlibrary{automata, positioning, arrows, calc}
\usepackage{subcaption}
\usepackage{wrapfig}

\newcommand{\vo}[1]{\boldsymbol{#1}}

\renewcommand{\b}{\vo{b}}

\newcommand{\elab}[1]{\label{eqn:#1}}
\newcommand{\eqn}[1]{(\ref{eqn:#1})}

\newcommand{\flab}[1]{\label{fig:#1}}
\newcommand{\fig}[1]{Fig.\ref{fig:#1}}

\newcommand{\etal}{\textit{et al.}}

\tikzstyle{block} = [draw, fill=white, rectangle, minimum height=3em, minimum width=4em, thick]
\tikzstyle{sum} = [draw, fill=white, circle, thick]
\tikzstyle{input} = [coordinate]
\tikzstyle{noise} = [coordinate]
\tikzstyle{output} = [coordinate]
\tikzstyle{output1} = [coordinate]
\tikzstyle{disturbance} = [coordinate]
\tikzstyle{pinstyle} = [pin edge={to-,thin,black}]

\title{\LARGE \bf
$\mathcal{H}_\infty$ Optimal Navigation in the Cislunar Space with LFT Models
}

\author{Tanay Kumar$^{1}$ and Raktim Bhattacharya$^{2}$
\thanks{This work was supported by AFOSR Grant FA9550-22-1-0539.}
\thanks{$^{1}$Graduate Student, Aerospace Engineering, Texas A\&M University. Email: ktanay@tamu.edu.}
\thanks{$^{2}$Professor, Aerospace Engineering, Texas A\&M University. Email: raktim@tamu.edu.}
}

\begin{document}

\maketitle
\thispagestyle{empty}
\pagestyle{empty}

\begin{abstract}
Navigation in the cislunar domain presents significant challenges due to chaotic and unmodeled dynamics, as well as state-dependent sensor errors. This paper develops a robust estimation framework based on Linear Fractional Transformation (LFT) models, and state estimation in $\mathcal{H}_\infty$ and $\mu$ synthesis framework to address these challenges. The cislunar dynamics are embedded into an LFT form that captures nonlinearities in the gravitational model and state-dependent sensor errors as structured uncertainty. A nonlinear estimator is then synthesized in the $\mathcal{H}_\infty$ sense to ensure robust performance guarantees in the presence of the stated uncertainties. Simulation results demonstrate the effectiveness of the estimator for navigation in a surveillance constellation. 
\end{abstract}

\section{Introduction}
\subsection{Motivation}
The Cislunar region has become increasingly important due to renewed interest in sustained lunar exploration, commercial activity, and strategic infrastructure development. As spacecraft begin to operate beyond Earth orbit with limited ground-based oversight, reliable onboard navigation becomes essential for mission autonomy and safety. High-confidence navigation is especially critical in the cislunar environment, where weak multi-body gravitational fields and limited sensing opportunities amplify the impact of model and sensor uncertainties. 

\subsection{Current Approaches \& Limitations}
Current approaches to cislunar navigation primarily rely on variants of the Extended Kalman Filter or Unscented Kalman Filter, often combined with optical measurements and onboard propagators. Although effective in many scenarios, these probabilistic estimation techniques provide limited guarantees under worst-case uncertainty. Kalman-based filters such as the Extended Kalman Filter (EKF) and Unscented Kalman Filter (UKF) assume accurate prior knowledge of the process and measurement models, along with Gaussian noise statistics. In practice, however, the cislunar environment is characterized by highly nonlinear dynamics, partial observability, and sensor degradations—conditions under which these assumptions often break down. As a result, the performance of such filters can degrade unpredictably in the presence of large modeling errors, abrupt disturbances, or non-Gaussian noise. A recent NASA report \cite{carpenter2025navigation} highlights that while the EKF has been a foundational tool in space navigation, it is based on several ad hoc assumptions that can lead to misuses and misunderstandings. Specifically, the report notes that the EKF’s reliance on local linearization and Gaussian noise assumptions can result in filter divergence or suboptimal performance in the presence of significant nonlinearities or unmodeled dynamics. 

\subsection{A Case for $\mathcal{H}_\infty$ and $\mu$-synthesis Framework for Cislunar Navigation} 
These challenges motivate the need for robust estimation frameworks, such as those based on $\mathcal{H}_\infty$ \cite{zhou1998essentials} and $\mu$  synthesis \cite{packard1988s}, which can provide guaranteed performance bounds based on norm-bounded uncertainty representation. Norm-bounded uncertainty models are particularly well-suited for safety-critical systems, such as spacecraft navigation, where worst-case guarantees are more critical than average-case performance. Unlike probabilistic models, which rely on well-characterized noise statistics, norm-bounded approaches remain effective when such information is unavailable or unreliable. They also allow structured representations of uncertainty (e.g., parameter drift or state-dependent sensor errors) that are difficult to model probabilistically. Moreover, norm-bounded frameworks align naturally with formal verification and validation requirements by ensuring performance guarantees over all admissible uncertainties, making them particularly valuable for certifiable autonomous systems operating in the cislunar environment.


Despite their advantages in handling worst-case uncertainties, $\mathcal{H}_\infty$ and $\mu$ analysis frameworks exhibit several limitations. First, $\mathcal{H}_\infty$ methods are inherently conservative, as they optimize performance against the worst possible disturbance consistent with the uncertainty bounds, which can potentially lead to overdesign in practical scenarios. Additionally, these methods typically assume norm-bounded, time-invariant uncertainties, which may not fully capture the complexity of time-varying or structured uncertainties encountered in cislunar navigation. $\mu$ analysis partially alleviates this conservatism by exploiting the structure of uncertainties. However, its computational complexity increases rapidly with the number of uncertainty blocks, thereby limiting its scalability. Moreover, both frameworks rely on linear fractional representations, which may introduce approximation errors when embedding strongly nonlinear dynamics. Finally, in general, the synthesis procedures often result in high-order filters or controllers that can be difficult to implement in resource-constrained systems. 

However, a key advantage of the  $\mathcal{H}_\infty$ and $\mu$ framework lies in the precise embedding of the Circular Restricted Three-Body Problem (CR3BP) \cite{koon2000dynamical} dynamics into a Linear Fractional Transformation (LFT) structure \cite{zhou1998essentials}, eliminating the need for local linearization or approximation. The CR3BP equations are algebraic-rational in form, involving nonlinearities that depend on the inverse of the distance to the primaries (e.g., $1/r_{13}$, $1/r_{23}$), the primaries being the distances of the space object (the third body) from the Earth and the Moon. These rational expressions can be modeled exactly using static nonlinear operators, which are treated as uncertainty blocks within the LFT formalism. By encapsulating these nonlinear terms as structured, bounded uncertainties around a linear time-invariant core, the LFT representation retains the full fidelity of the nonlinear dynamics while enabling robust synthesis via $\mathcal{H}_\infty$ or $\mu$ tools. This capability makes the LFT framework particularly well-suited for cislunar navigation, where high-fidelity modeling is essential and small-angle or local-linear approximations are insufficient. Also, we utilize the framework to design an observer that relies on the design of a gain matrix with fixed dimensions, thereby eliminating concerns associated with the numerical integration of higher-order dynamical systems in onboard computers.

The use of LFT modeling combined with $\mathcal{H}_\infty$ and $\mu$ synthesis tools is motivated by the growing need for certifiable robustness in autonomous, safety-critical space systems. As these systems increasingly operate without human oversight, design methodologies rooted in operator-theoretic analysis and structured uncertainty, such as those enabled by the LFT framework, offer a principled path toward establishing trust in navigation and control software under verification and validation constraints.  

LFT-based modeling has already seen wide adoption in the verification and validation of aerospace systems. For instance, in the development of flight control systems for re-entry vehicles, LFT models have been employed to capture nonlinear dynamics and parametric uncertainties, enabling robust control synthesis and certification under stringent safety requirements. Similarly, NASA has utilized LFT formulations in the robustness analysis and reliable flight regime estimation of aircraft, enabling engineers to assess closed-loop stability margins under varying aerodynamic conditions and actuator saturations. The European Space Agency (ESA) has also integrated LFT modeling into its robust design and validation workflows, particularly for spacecraft attitude control systems that are subject to high-fidelity parameter variations. These applications highlight the practical relevance of LFT theory in developing certifiable, high-confidence designs for safety and mission-critical aerospace systems \cite{shin2008robustness, di2010integrated, roux2012robust}. 

\subsection{Contributions of the Paper}
Despite the maturity of LFT-based robust control techniques in aerospace applications, their application to spacecraft navigation remains limited, presenting an opportunity to enhance the reliability and certifiability of autonomous navigation systems. This paper demonstrates how such tools can be used in spacecraft navigation in cislunar space.

Specifically, the following are the four main contributions of this paper.
\begin{enumerate} 
\item It develops an exact LFT representation of the Circular Restricted Three-Body Problem (CR3BP), preserving nonlinear gravitational dynamics and sensing without approximation. 
\item It models state-dependent sensing noise, a common phenomenon in cislunar navigation scenarios, within the LFT framework.
\item It synthesizes a nonlinear state estimator using $\mathcal{H}_\infty$ and $\mu$ synthesis tools, offering formal robustness guarantees under both structured and unstructured uncertainties. 
\item Finally, the paper demonstrates the effectiveness of this approach through simulation studies of spacecraft navigation in cislunar space, highlighting its potential for certifiable autonomous navigation systems.
\end{enumerate}

This paper presents a framework based on LFT for modeling nonlinear systems, facilitating the development of robust observers through $\mathcal{H}_\infty$ and $\mu$ methods. Unlike our previous research \cite{nychka2024optimal}, which used polytopic uncertainty and function substitution to derive affinely parameter-dependent linear models, leading to a more complex polytopic uncertainty description, this approach formulates the CR3BP dynamics and sensor model using an exact LFT representation. This method enables the straightforward application of $\mathcal{H}_\infty$ and $\mu$ synthesis techniques for designing nonlinear observers.

\section{LFT Modeling of Dynamics and Sensing}
\subsection{The Circular Restricted Three-Body Dynamics}
We consider the classical circular restricted three-body problem (CR3BP), which models the motion of a small body under the gravitational influence of two larger bodies in circular orbits around their common center of mass. The two larger bodies, known as the primary bodies, are significantly more massive than the small body and are assumed to follow fixed circular orbits due to their mutual gravitational attraction. The small body, which is considered to have negligible mass, exerts no influence on the motion of the primary bodies. In the cislunar context, as shown in \fig{cislunar}, the Earth and the Moon serve as the primary bodies, with the small body representing a spacecraft or another resident space object.

\begin{figure}[ht!]
\includegraphics[width=0.5\textwidth]{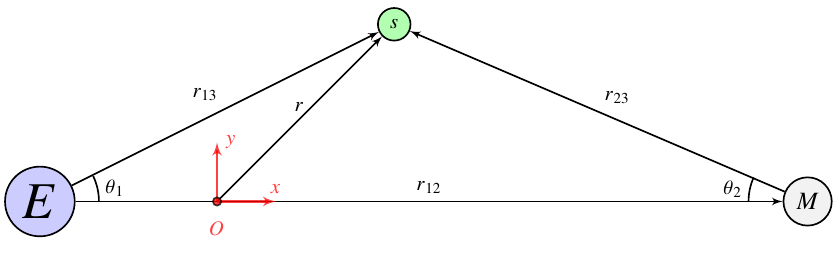}
\caption{Cislunar Rotating Reference Frame, with $E$ (Earth), $M$ (Moon), and $s$ (Satellite).}
\flab{cislunar}
\end{figure}

The nonlinear equations of motion in \textit{two dimensional} space are given by
\begin{subequations}
\begin{align}
\ddot{x} - 2\Omega\dot{y} - \Omega^2x &= -\frac{\mu_1}{r_{13}^3}(x + \pi_2r_{12}) - \frac{\mu_2}{r_{23}}(x - \pi_1r_{12}), \\ 
\ddot{y} - 2\Omega\dot{x} - \Omega^2y &= -\frac{\mu_1}{r_{13}^3}y -\frac{\mu_2}{r_{23}^3}y, 
\end{align}
\elab{dynamics}
\end{subequations}
with constants \(\pi_1\), \(\pi_2\), \(\mu\), \(\mu_1\), \(\mu_2\), \(\Omega\), and \(T\) given by,
\begin{equation}\left.
\begin{split}
&\pi_1 = \frac{m_1}{m_1 + m_2}, \quad \pi_2 = \frac{m_2}{m_1 + m_2}, \quad \mu_1 = Gm_1, \\ 
&\quad \mu_2 = Gm_2,\quad \mu = \mu_1+\mu_2,\quad \Omega = \sqrt{\frac{\mu}{{r_{12}}^3}}, \quad T = \frac{2\pi}{\Omega},
\end{split}\right\}
\elab{constants}
\end{equation}
where $m_1$ and $m_2$ are the masses of the primaries bodies, $G$ is the universal Gravitational constant, and $r_{ij}$ is the distance between bodies $i$ and $j$. The third dimension is excluded for simplicity and does not limit the generality of the result.

For numerical stability, the equations are expressed in a nondimensionalized form by introducing
\begin{equation}
\begin{split}
\boldsymbol{\lambda} = \frac{\mathbf{r}}{r_1{}_2} &= \Bar{x}\hat{i} + \Bar{y}\hat{j}, \\
\boldsymbol{\sigma} = \frac{\mathbf{r_{13}}}{r_1{}_2} &= (\Bar{x} + \pi_2)\hat{i} + \Bar{y}\hat{j}, \\
\boldsymbol{\psi} = \frac{\mathbf{r_{23}}}{r_1{}_2} &= (\Bar{x} + \pi_2 - 1)\hat{i} + \Bar{y}\hat{j},
\end{split}
\elab{vectors}
\end{equation}
where the normalized positions are given by $\Bar{x} := \frac{x}{r_1{}_2}$ and $\Bar{y} := \frac{y}{r_1{}_2}$ respectively, and the time is non-dimensionalized with respect to $t_c := \sqrt{\frac{r_1{}_2^3}{\mu}}$.

Substituting these in \eqn{dynamics}, we get the following nondimensionalized equations of motion:
\begin{subequations}
\begin{align}
\ddot{\Bar{x}} - 2\dot{\Bar{y}} - \Bar{x} &= -\frac{1 - \pi_2}{\sigma^3}(\Bar{x} + \pi_2) - \frac{\pi_2}{\psi^3}(\Bar{x} - 1 + \pi_2), \\
\ddot{\Bar{y}} + 2\dot{\Bar{x}} - \Bar{y} &= -\frac{1 - \pi_2}{\sigma^3}\Bar{y} - \frac{\pi_2}{\psi^3}\Bar{y},
\end{align}
\elab{normalized_dynamics}
\end{subequations}
where \(\sigma\) is the \(|\boldsymbol{\sigma}|\) and  \(\psi\) is the \(|\boldsymbol{\psi}|\). 

The nonlinearities in \eqn{normalized_dynamics} are in $\sigma(\bar{x},\bar{y})$ and  $\psi(\bar{x},\bar{y})$, which can be expressed as a parameter dependent system 

\begin{align}
\notag \frac{d}{dt}\begin{pmatrix}\bar{x} \\ \bar{y} \\ \dot{\bar{x}} \\ \dot{\bar{y}} \end{pmatrix} &=
\begin{bmatrix}
0 & 0&  1& 0\\
0 & 0&  0& 1\\
a_{31}(\sigma,\psi) & 0&  0& 2\\
0 & a_{42}(\sigma,\psi) & -2 & 0
\end{bmatrix}\begin{pmatrix}\bar{x} \\ \bar{y} \\ \dot{\bar{x}} \\ \dot{\bar{y}} \end{pmatrix} \\  &+\begin{pmatrix} 0 \\ 0 \\b_3(\sigma,\psi)\\0
\end{pmatrix},
\elab{lft_dynamics}
\end{align}
where
\begin{align*}
a_{31}(\sigma,\psi) &:= (\pi_2 - 1)/\sigma^3 - \pi_2/\psi^3 + 1,\\
a_{42}(\sigma,\psi) &:= (\pi_2 - 1)/\sigma^3 - \pi_2/\psi^3 + 1,\\
b_3(\sigma,\psi) &:= \pi_2(1-\pi_2)(1/\psi^3-1/\sigma^3).
\end{align*}

Assuming $\sigma$ and $\psi$ are bounded (can be verified using the Jacobi Constant or Poincar\'e maps), i.e., $\sigma\in[\sigma_\text{min},\sigma_\text{max}]$, and $\psi\in[\psi_\text{min},\psi_\text{max}]$, we can treat them as multiplicative structured uncertainty, i.e.,
$$
\sigma = \bar{\sigma}(1+\delta_\sigma\tilde{\sigma}), \text{ and } \psi = \bar{\psi}(1+\delta_\psi\tilde{\psi}),
$$
where $\delta_\sigma \in [-1,1]$ and $\delta_\psi \in [-1,1]$, and $\bar{\sigma}$ and $\bar{\psi}$ are nominal values, often taken as the average of the extreme values. We can express the nonlinear uncertain `parameters' in LFT form, resulting in the following uncertain linear dynamical system (without any approximations). This can be accomplished using MATLAB's Robust Control Toolbox~\cite{matlabRobustControlToolbox}, which provides functions for LFT modeling and synthesis.

\subsection{The Sensing Model and Associated Uncertainty}
Motivated by the Orion optical navigation (OpNav) system on Artemis 1 \cite{inman2024artemis, inman2024demonstration, christian2012onboard}, we consider designing a celestial navigation system using only bearing measurements. Specifically, we assume angles $\theta_1$ and $\theta_2$, as shown in \fig{cislunar}, are available for navigation. We incorporate bearing measurements as sines and cosines of measured bearing angles to enable an LFT modeling of sensors. Therefore, the equations for the four measurements $y_{m_1},\cdots,y_{m_4}$ are
\begin{equation}
\begin{split}
y_{m_1} &= \sin{\theta_1} = \frac{\Bar{y}}{\sigma}, \\
y_{m_2} &= \cos{\theta_1} = \frac{\Bar{x} + \pi_2}{\sigma}, \\
y_{m_3} &= \sin{\theta_2} = \frac{\Bar{y}}{\psi},\\
y_{m_4} &= \cos{\theta_2} = \frac{\Bar{x} + \pi_2 - 1}{\psi},
\end{split}
\elab{sensor}
\end{equation}
which can be compactly represented as
\begin{equation}
y_m = \begin{bmatrix} 
    0&    1/\sigma & 0 & 0\\
    1/\sigma&   0&   0& 0\\
    0& 1/\psi& 0& 0\\
  -1/\psi&     0&  0& 0
\end{bmatrix} \begin{pmatrix}\bar{x} \\ \bar{y} \\ \dot{\bar{x}} \\ \dot{\bar{y}} \end{pmatrix}  +  \begin{pmatrix} 0\\ \frac{\pi_2}{\sigma}\\0\\-\frac{(\pi_2 - 1)}{\psi}\end{pmatrix} + v,
\elab{sensor_matrix}
\end{equation}
where $v$ is the measurement noise vector. The sensor model in \eqn{sensor} is also nonlinear in $\sigma$ and $\psi$, which can be expressed in LFT form.

In optical navigation systems, sensor measurement noise $v$ is strongly dependent on the range of the observed object due to the degradation of the signal-to-noise ratio as the distance increases. As the target's apparent brightness decreases, the precision of image-based measurements, such as centroids or limb features, deteriorates. This effect is especially pronounced in the cislunar environment, where the large and varying distances between Earth and the Moon amplify the impact of photon-limited sensing. To model this effect within a robust estimation framework, the sensor noise is expressed as a range-weighted noise. Consequently, the measurement noise in \eqn{sensor_matrix} can be expressed as
\begin{equation}
v =  \textbf{blkdiag}\begin{pmatrix}W_1(\sigma), & W_2(\sigma), & W_3(\psi), & W_4(\psi)\end{pmatrix}\bar{v},    
\elab{noise_model}
\end{equation}
where $\bar{v}$ is a unit-norm exogenous noise, and \( W_k(\cdot) \) is a range-dependent weighting function that captures the growth of sensor noise with distance. A common and physically motivated choice is \( W_k(r) = \sqrt{\alpha}_k r_k \), so that the noise energy scales quadratically with range: $\|(W_k(r_k)\bar{v_k}\|^2 = \alpha_k r_k^2$, where $r_k$ is either $\sigma$ or $\psi$, depending on the sensor. This formulation preserves the structure necessary for robust $\mathcal{H}_\infty$ and $\mu$ synthesis while accurately capturing the degradation in measurement precision at large distances. Such a representation is supported by optical navigation data from missions like Artemis I, where increased range to Earth or Moon led to visibly higher measurement uncertainty \cite{inman2024artemis}. Incorporating this range-weighted noise model into the estimator design provides performance guarantees that consider the state-dependent sensing accuracy in cislunar operations.

The nonlinear measurement model in \eqn{sensor}, along with the state-dependent noise modeled in \eqn{noise_model} can also be expressed in LFT form using MATLAB's Robust Control Toolbox~\cite{matlabRobustControlToolbox}.

\section{Robust State-Estimation Problem as an LFT Inerconnection}
The dynamics in \eqn{lft_dynamics} and the sensor model in \eqn{sensor_matrix} can be expressed in a linear parameter varying (LPV) form, where the parameters are the nonlinear terms $\sigma$ and $\psi$. The state-space representation of the system can be written as

\begin{equation}
\begin{split}
\dot{x} &= A(\rho)x + B_ww + b(\rho),\\
y_m &= C_y(\rho)x + D_w(\rho)w + d(\rho),\\
z &= C_zx,
\end{split}
\elab{lft_system}
\end{equation}
where \(x\) is the state vector, \(y_m\) is the measurement vector, \(w\) is the exogenous inputs including disturbances and sensor noise. The vector $z$ represents the output of interest, which for navigation purposes is the position $(x,y)$. The vector $\rho$ represents the parameters $\sigma$ and $\psi$, i.e., $\rho := \begin{pmatrix}\sigma & \psi\end{pmatrix}^T$. The system matrices \(A(\rho)\), \(b(\rho)\), \(C_y(\rho)\), \(D_w(\rho)\), $d(\rho)$ are nonlinear functions of $\rho$. The nonlinearities in the system can be expressed as structured uncertainties within the LFT framework, enabling the application of robust estimation techniques.

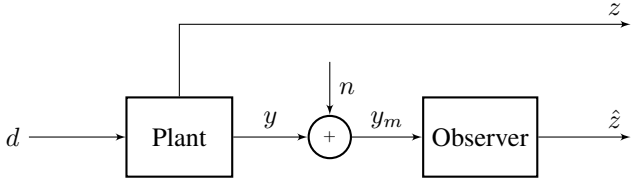
\begin{figure}

\centering

\begin{tikzpicture}[auto, node distance=2cm, >=latex']
    \node [block, right of=d, node distance=2cm] (P) at (2,-0.25){Plant};
    \node [sum, right of=P, node distance=2cm, thick] (sum1) {\tiny $+$};
    \node [noise, above of=sum1, node distance=1cm] (noise) {};
    \node [block, right of=sum1, node distance=2cm] (observer) {Observer};
    \node [output, right of=observer] (zh) {};
    \node [output1, above of=zh, node distance=1.5cm] (z) {};
    
    \node [input, left of=P, yshift=0mm] (d) {$d$};
    \node [input, left of=P, yshift=-2mm] (u) {$u$};
    \coordinate [left of=P, yshift=-2mm, xshift=0.6cm] (u1) {};
    \coordinate [below of=u1, node distance = 10mm] (u2) {};
    
    \draw [->] (P) -- node {$y$} (sum1);
    \draw [->] (noise) -- node {$n$} (sum1);
    \draw [->] (sum1) -- node {$y_m$} (observer);
    \draw [->] (observer) -- (zh) node [at end, above, above left] {$\hat{z}$};
    \draw [->] (P) |- (z) node [at end, above left] {$z$};
    \draw [->] (d) node[left] {$d$} -- ++(1.29,0);

\end{tikzpicture}
\caption{System interconnection for designing and implementing proposed estimators.}
\flab{Block Diagram}
\end{figure}

The objective is to design a state estimator that provides accurate estimates of the outputs of interest \(z\) based on the measurements \(y_m\), while ensuring robustness against uncertainties in the parameters \(\rho\) and exogenous inputs \(w\). The schematic of the estimator is shown in \fig{Block Diagram}. The estimator can be designed using the $\mathcal{H}_\infty$ and $\mu$ synthesis framework, which enables the incorporation of structured uncertainties into the system.

Following the derivation of standard full-order $\mathcal{H}_{\infty}$ observer \cite{duan2013lmis} (pg. 293), we propose the following observer structure for the system in \eqn{lft_system},
\begin{equation}
\dot{\hat{x}} = \big(A(\rho) + LC_y(\rho)\big)\hat{x} - L\big(y_m-d(\rho)\big) + b(\rho),
\elab{obs}
\end{equation}
which is slightly different from the standard observer structure due to the presence of the terms $b(\rho)$ and $d(\rho)$ in \eqn{lft_system}. In LPV models, the parameter $\rho$ is assumed to be known at runtime -- either via direct measurement or estimation.

Defining the error as $e(t) := x(t) - \hat{x}(t)$, we get the following error dynamics
$$
\dot{e} = \big(A(\rho) + LC_y(\rho)\big)e + \big(B_w + LD_w(\rho)\big)w,
$$
which is in the form shown in Duan \etal \cite{duan2013lmis} (pg. 293). The estimated quantity of interest is $\hat{z} := C_z\hat{x}$, and the error in the estimate $\tilde{z}$ is related to the error $e$ as $\tilde{z} = C_ze$. Therefore, the dynamical system that relates exogenous input $w(t)$ to estimation error $\tilde{z}(t)$ is given by
\begin{equation}
\begin{split}
\dot{e} &= \big(A(\rho) + LC_y(\rho)\big)e + \big(B_w + LD_w(\rho)\big)w,\\
\tilde{z} &= C_ze.
\end{split}
\elab{estim-error}
\end{equation}
The objective is to determine $L$ such that $\|\tilde{z}(t)\|_2$ is minimized, which is achieved by minimizing the $\mathcal{H}_\infty$ norm of the system in \eqn{estim-error} \cite{doyle2013feedback}.

To synthesize the observer gain $L$, we can express the system in LFT form. The LFT representation of the system in \eqn{estim-error} can be written as
\begin{equation}
\begin{split}
  \dot{e} &= A(\rho)e + B_ww + \tilde{u}, \\
  \tilde{y}_m &= C_y(\rho)e + D_w(\rho) w,\\
  \tilde{u} &= L\tilde{y}_m,\\
  \tilde{z} &= C_ze.
\end{split}
\elab{lft_estim}
\end{equation}

The system in \eqn{lft_estim} can be expressed in the LFT form as shown in \fig{lft_estim}.
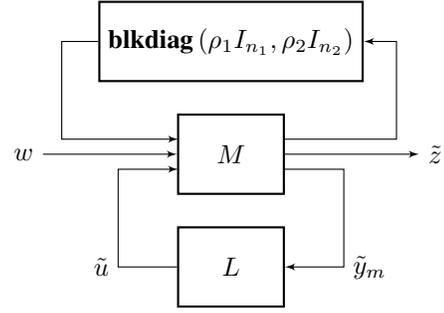
\begin{figure}
  \centering
\begin{tikzpicture}[auto, node distance=1.5cm, >=latex']
    
    \node [block] (P) at (0,0){$M$};
    \node [block, above of=P] (delta) {$\textbf{blkdiag}\left(\rho_1 I_{n_1},\rho_2 I_{n_2}\right)$};
    \node [block, below of=P] (L) {$L$};

    \node [input, left of=P, yshift=2mm] (w1) {};
    \node [input, left of=P, yshift=0mm] (w2) {};
    \node [input, left of=P, yshift=-2mm] (u) {};

    \node [output, right of=P, xshift=-.8cm, yshift=2mm] (z1) {};
    \node [output, right of=P, yshift=0mm] (z2) {};
    \node [output, right of=P, xshift=-.8cm, yshift=-2mm] (ym) {};

    \draw [->] (delta) -- ++(-2.25,0) |- (w1) -- ++(0.8,0);
    \draw [<-] (P) -- ++(-2.5,0) node[left] {$w$};
    \draw [->] (L) -| node[left] {$\tilde{u}$} (u) -- ++(0.8,0);

    \draw [->] (z1) -- ++(1.5,0) |- (delta);
    \draw [->] (ym) -- ++(0.8,0) |- node [midway,right] {$\tilde{y}_m$}(L);
    \draw [->] (P) -- ++(2.5,0) node[right] {$\tilde{z}$};

\end{tikzpicture}
\caption{LFT interconnection for designing robust $\mathcal{H}_\infty$ estimator.}
\flab{lft_estim}
\end{figure}

The observer $L$ can be designed by expressing $\rho$ as \texttt{ureal} objects in MATLAB's Robust Control Toolbox, which constructs the uncertainty blocks in the LFT representation in \fig{lft_estim}. The observer $L$ can be synthesized using the \texttt{hinfstruct($\cdots$)} for the $\mathcal{H}_\infty$ synthesis, and \texttt{musyn($\cdots$)} for the $\mu$ synthesis. The resulting observer gain $L$ is a fixed-dimension matrix.

\section{Simulation Results}
\subsection{Simulation Setup}
We consider a use case in which a satellite, positioned within a constellation, performs self-localization using bearing and range measurements. The simulations are performed with initial conditions,
\begin{align}
x(0) &:= \begin{pmatrix}0.87 & 0 &0 &-1.48270\end{pmatrix}^T,\\
\hat{x}(0) &:= \begin{pmatrix} 0.65 & -0.1 & -2 & -2.0 \end{pmatrix}^T.
\end{align}

Note a significant difference between the actual and estimated initial conditions. The limits for $\rho$ are obtained from the constellation trajectories and are listed in table \ref{tab:rho_bounds}.
\begin{table}
\begin{center}
\begin{tabular}{l||r|r|}
    & $\rho_1:=\sigma$ & $\rho_2:=\psi$ \\[1mm] \hline
Min & 0.1289 & 0.1218 \\
Max & 0.9005 & 1.9005\\
\end{tabular}
\caption{Bounds for parameters.}
\label{tab:rho_bounds}
\end{center}
\end{table}

We consider exogenous signal vector  $$w(t) = \begin{pmatrix} d_x(t) & d_y(t) & n_1(t) & n_2(t) & n_3(t) & n_4(t)\end{pmatrix}^T,$$ where $d_x(t)$ and $d_y(t)$ represents the perturbation accelerations acting on the spacecraft (caused by the Sun, Jupiter, Earth’s gravitational harmonics, etc.) in the $x$ and $y$ directions respectively, and $n_1(t)$, $n_2(t)$, $n_3(t)$, $n_4(t)$ are the sensor noises in the measurements. Therefore, 
\begin{align*}
B_w(\rho) &:= \begin{bmatrix}
0_{2\times 2} & 0_{2\times 4}\\ I_2 & 0_{2\times 4}
\end{bmatrix},\\ 
D_w(\rho) &:= \begin{bmatrix}
0_{4\times 2} & \textbf{blkdiag}\left(W_1(\sigma)I_2,W_2(\psi)I_2\right)
\end{bmatrix},
\end{align*}
where 
\begin{equation}
\begin{split}
W_1(\sigma) &:= \eta_\text{min} + \frac{\sigma - \sigma_\text{min}}{\sigma_\text{max} - \sigma_\text{min}}(\eta_\text{max} - \eta_\text{min}),\\
W_2(\psi) &:= \eta_\text{min} + \frac{\psi - \psi_\text{min}}{\psi_\text{max} - \psi_\text{min}}(\eta_\text{max} - \eta_\text{min}),
\end{split}
\elab{sensor-noises}
\end{equation}
where $\eta_\text{min}$ and $\eta_\text{max}$ represent the minimum and maximum noises that vary linearly with the range. In the simulations, we consider $\eta_\text{max} = 500$ arcsec and $\eta_\text{min} = 50$ arcsec \cite{wu2022autonomous}. We also need the range measurements to determine $\rho$ in real-time -- to implement the estimator in \eqn{obs} -- and consider a similar sensor noise model as in \eqn{sensor-noises} with minimum and maximum errors equal to $400$ and $4000$ km, respectively. Sensor noises are generated as uniform random numbers within the specified intervals. We investigate the performance of the estimators with uniform white noise and band-limited uniform white noise (1 rad/s). 

Process noise is considered uniform white noise within the $[-0.01, 0.01]$ range. This range is relatively large in the context of normalized CR3BP dynamics and is selected to assess the robustness of the observers without any specific physical justification.

Since we are interested in estimating only the spacecraft's position in this example, $C_z := \begin{bmatrix}I_2 & 0_{2\times 2}\end{bmatrix}$. 
The rest of the matrices in \eqn{lft_estim} are defined in \eqn{lft_dynamics} and \eqn{sensor_matrix}, and $I_n$ is the $n\times n$ identity matrix.

\subsection{Results and Discussion}

\begin{figure}[h!]\centering
  \includegraphics[width=0.4\textwidth]{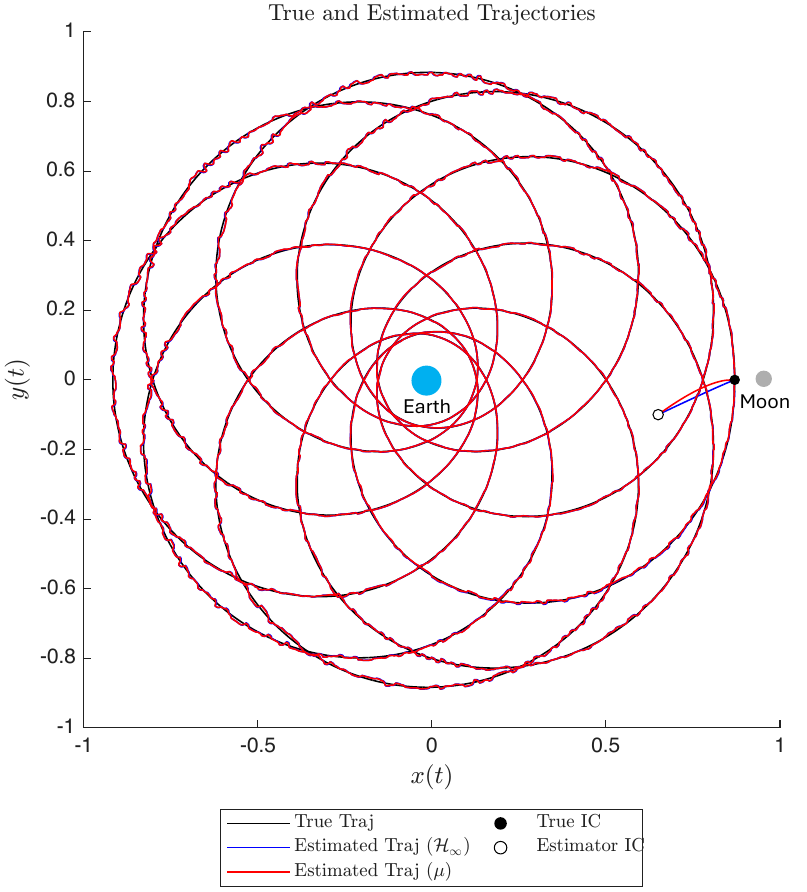}
  \caption{Estimated trajectories of the spacecraft using $\mathcal{H}_\infty$ and $\mu$ observers.}
  \flab{estimated_trajectories}
\end{figure}

\begin{figure}[h!] \centering
  \includegraphics[width=0.4\textwidth]{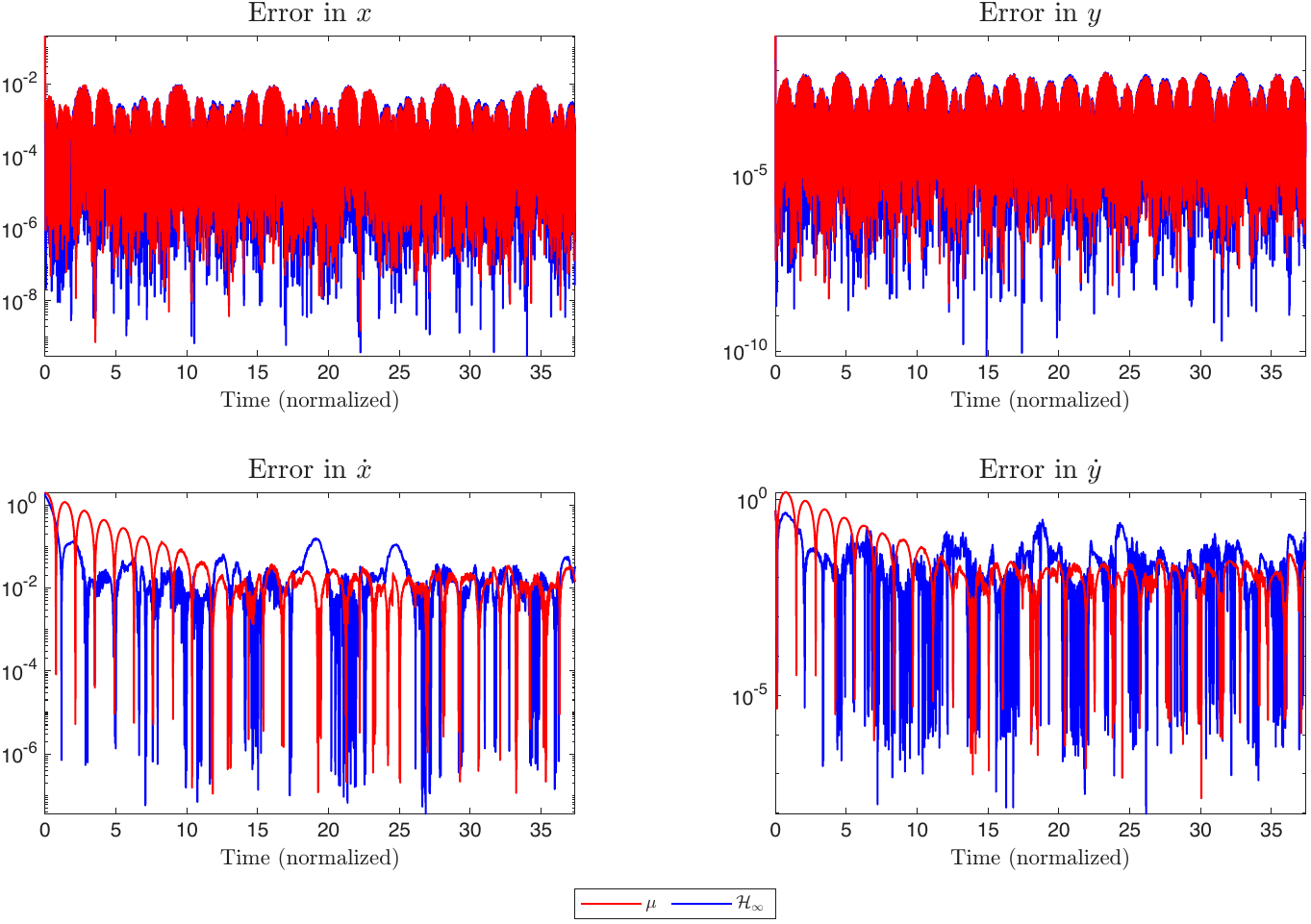}
  \caption{Estimation errors of the spacecraft position using $\mathcal{H}_\infty$ and $\mu$ observers. Errors are in normalized state and time.}
  \flab{estimation_errors}
\end{figure}
  
The estimated trajectories of the spacecraft using the $\mathcal{H}_\infty$ and $\mu$ observers are shown in \fig{estimated_trajectories}. In this simulation, the true trajectory corresponds to a surveillance orbit derived from an Earth-Moon resonant orbit, which is commonly used for persistent monitoring in the cislunar region. Such orbits are designed to periodically revisit specific locations relative to the Earth and Moon, providing favorable trajectories for cislunar surveillance \cite{frueh2021cislunar, gupta2023constellation}. The robust observers are tasked with localizing the spacecraft as it follows the surveillance orbit, using only bearing and range measurements despite the presence of significant nonlinearities and state-dependent sensor noise. The results demonstrate that both observers are able to accurately reconstruct the surveillance orbit, maintaining close agreement with the true trajectory throughout the simulation.

The estimation errors are shown in \fig{estimation_errors}. The results indicate that both observers provide accurate estimates of the spacecraft's position, with the $\mu$ observer exhibiting slightly better long-term performance in terms of estimation accuracy. The errors remain bounded, demonstrating the robustness of the proposed observers against uncertainties in the system dynamics and sensor measurements.

As evident from \fig{estimated_trajectories} and \fig{estimation_errors}, the state-dependent sensor-noise model leads to larger estimation errors when the satellite is farther from the Earth and the Moon. This is a direct consequence of the range-weighted noise formulation, where the measurement noise increases with the distance to the observed bodies. When the spacecraft is at greater ranges, the signal-to-noise ratio of the optical navigation measurements degrades, resulting in less precise bearing and range information. Consequently, the observer's ability to accurately estimate the states is diminished during these periods, leading to increased estimation errors. This effect is particularly pronounced in the cislunar environment, where the distances can vary significantly throughout an orbit. The simulation results confirm that the robust observers maintain bounded errors, but the magnitude of these errors is correlated with the satellite's distance from the primary bodies, reflecting the physical limitations imposed by state-dependent sensor noise. We observe that the $\mathcal{H}_\infty$ observer exhibits slightly larger errors compared to the $\mu$ observer during the periods when the spacecraft is farther from the Earth and Moon. This difference can be attributed to the $\mu$ observer's ability to handle structured uncertainties more effectively, leading to improved robustness in the presence of state-dependent sensor noise. The transient errors are, however, better for the $\mathcal{H}_\infty$ observer.

\section{Conclusions}
This paper presents a robust state estimation framework for spacecraft navigation in cislunar space using Linear Fractional Transformation (LFT) modeling and $\mathcal{H}_\infty$ and $\mu$ synthesis techniques. The proposed approach captures the nonlinear dynamics of the Circular Restricted Three-Body Problem (CR3BP) and models state-dependent sensor noise, enabling accurate and robust estimation of spacecraft states. 

The LFT-based robust estimation framework offers potentially superior performance for cislunar navigation by directly modeling nonlinearities and structured uncertainties without relying on local linearization or Gaussian noise assumptions. This approach provides worst-case performance guarantees through $\mathcal{H}_\infty$ and $\mu$ synthesis, handles state-dependent and structured uncertainties more naturally, and aligns with formal verification requirements for safety-critical systems. As a result, LFT-based methods deliver greater reliability and certifiability than traditional Kalman filters, especially in the highly uncertain and nonlinear cislunar environment. The limitation in this approach is the assumption of continuous measurements, which may not hold in practice due to sensor limitations or occlusions. This is a subject of our current research.

Kalman filtering remains valuable when sensor measurements are intermittent, as its predict-update structure allows state propagation during measurement gaps. This makes it suitable for cislunar missions, which are characterized by sensing occlusions. However, Kalman filters (and its nonlinear/non-Gaussian extensions) lack formal robustness guarantees under nonlinear dynamics or non-Gaussian noise. In contrast, LFT-based $\mathcal{H}_\infty$ and $\mu$ observers assume continuous measurements but offer guaranteed robustness to nonlinearities and structured uncertainties. Extensions using event-triggered control can address intermittent measurement scenarios~\cite{Heemels2021}, and is the focus of our future work.

\section*{Acknowledgments}
This work is supported by AFOSR grant FA9550-22-1-0539 with Dr. Erik Blasch as the program director.

\bibliographystyle{unsrt}
\bibliography{references}

\end{document}